\documentclass[11pt]{amsart}
\usepackage{amssymb, amscd}
\usepackage{latexsym,epsfig}
\usepackage{pst-all}
\numberwithin{equation}{section}

\newtheorem{theorem}{Theorem}[section]
\newtheorem{lemma}[theorem]{Lemma}
\newtheorem{proposition}[theorem]{Proposition}
\newtheorem{corollary}[theorem]{Corollary}

\newtheorem{definition-lemma}[theorem]{Definition-Lemma}

\theoremstyle{definition}
\newtheorem{definition}[theorem]{Definition}
\newtheorem{example}[theorem]{Example}

\theoremstyle{remark}

\begin{document}

\title{Prym varieties associated to graphs}

\author{Rudi Salomon}
\address{Faculteit Wiskunde en Informatica, University of
Amsterdam, Plantage Muidergracht 24, 1018 TV Amsterdam, The
Netherlands.}
\email{rsalomon@science.uva.nl}


\begin{abstract} We present a Prym construction which associates
abelian varieties to vertex-transitive strongly regular graphs. As
an application we construct Prym-Tyurin varieties of arbitrary
exponent $\geq 3$, generalizing a result by Lange, Recillas and
Rochas.
\end{abstract}

\maketitle

\begin{section}{Introduction}

We describe a Prym construction which associates abelian varieties
to certain graphs. More precisely, given the adjacency matrix
$A=(a_{ij})_{i,j=1}^d$ of a vertex-transitive strongly regular
graph $\mathcal{G}$ along with a covering of curves
$p:C\rightarrow \mathbb{P}^1$ of degree $d$ and a labelling
$\{x_1,\ldots,x_d\}$ of an unramified fiber such that the induced
monodromy group of $p$ is represented as a subgroup of the
automorphism group of $\mathcal{G}$, we construct a symmetric
divisor correspondence $D$ on $C\times C$ which then serves to
define complementary subvarieties $P_{\!{}_+}$ and $P_{\!{}_-}$ of
the Jacobian $J(C)$. The correspondence $D$ is defined in such a
way that the point $(x_i,x_j)$ appears in $D$ with
multiplicity~$a_{ij}$, analogous to Kanev's construction in
\cite{Ka2}. The varieties $P_{\!{}_\pm}$ are given by
$P_{\!{}_\pm}={\rm ker}(\gamma-r_{\!{}_\mp}{\rm
id}_{J(C)})_{\!{}_0}$, where $r_{\!{}_\pm}$ are special
eigenvalues of $A$ and $\gamma$ is the endomorphism on $J(C)$
canonically associated to $D$ (i.e., sending the divisor class
$[x-x_0]$ to the class $[D(x)-D(x_0)]$). It is easy to show that
$$
(\gamma-r_{\!{}_+}{\rm id}_{J(C)})(\gamma-r_{\!{}_-}{\rm
id}_{J(C)})=0
$$ and $P_{\!{}_\pm}={\rm
im}(\gamma-r_{\!{}_\pm}{\rm id}_{J(C)})$. In particular, if $D$ is
fixed point free and $r_{\!{}_+}=1$, then $P_{\!{}_+}$ is a
Prym-Tyurin variety of exponent $1-r_{\!{}_-}$ for $C$. Given the
ramification of $p$ it is not hard to compute the dimension of
$P_{\!{}_\pm}$.

For a thorough definition of $D$ we consider the Galois closure
$\pi:X\rightarrow \mathbb{P}^1$ of $p$ and use the induced
representation ${\rm Gal}(\pi)\rightarrow {\rm Aut}(\mathcal{G})$
to construct symmetric correspondences $D_{\!{}_+}$ and
$D_{\!{}_-}$ on $X\times X$ (much the way M\'{e}rindol does in
\cite{Me}). With $C$ being a quotient curve of $X$, the
correspondence $D$ is derived from $D_{\!{}_\pm}$ taking quotients
and adding $r_{\!{}_\pm}\Delta_C$, where $\Delta_C$ is the
diagonal of $C\times C$; see section \ref{nongalois}. Given the
endomorphisms $\gamma_{\!{}_\pm}$ on $J(X)$ canonically associated
to $D_{\!{}_\pm}$, we show that ${\rm im}\,\gamma_{\!{}_\pm}$ and
$P_{\!{}_\pm}$ are isogenous.\\

The lattice graphs $L_2(n)$, $n\geq 3$ and their complements
$\overline{L_2(n)}$ offer important examples. For instance,
applying the method to $\overline{L_2(n)}$ and appropriate
coverings $C\rightarrow \mathbb{P}^1$ of degree $n^2$ with branch
loci of cardinality $2(l+2n-2)$ for $l\geq 1$, we obtain
$l$-dimensional Prym-Tyurin varieties of exponent $n$ for the
curves $C$; see section \ref{exponent}. We give a characterization
of these varieties and show that for $n=3$ they coincide with the
non-trivial Prym-Tyurin varieties of exponent 3 described by
Lange, Recillas and Rochas in \cite{LRR}.\\

\noindent \textbf{Conventions and notations.} The ground field is
assumed to be the field $\mathbb{C}$ of complex numbers. By a
covering of curves we mean a non-constant morphism of irreducible
smooth projective curves. The symbol $S_n$ denotes the symmetric
group acting on the letters $1,\ldots,n$ for $n\in \mathbb{N}$.

\end{section}

\begin{section}{Strongly regular graphs and matrices}
\label{matrices}

We start our discussion with the definition of a strongly regular
graph and its adjacency matrix and collect some properties of
such graphs. For additional information we refer to \cite{LW}.\\

By definition the set of {\it strongly regular} graphs ${\rm
SRG}(d,k,\lambda,\mu)$, $k> 0$, consists of the graphs
$\mathcal{G}$ with vertex set $\{v_1,\ldots,v_d\}$ such that
\begin{itemize}
\item[a)] the set $\Gamma(v_i)$ of vertices adjacent to $v_i$ has
exactly $k$ elements and $v_i\notin \Gamma(v_i)$;
\item[b)] for
any two adjacent vertices $v_i$,$v_j$ there are exactly $\lambda$
vertices adjacent to both $v_i$ and $v_j$;
\item[c)] for any two
distinct non-adjacent vertices $v_i$,$v_j$ there are exactly $\mu$
vertices adjacent to both $v_i$ and $v_j$.
\end{itemize} Let $A=(a_{ij})\in \{0,1\}^{d\times d}$
be the adjacency matrix of such a strongly regular graph
$\mathcal{G}$, i.e., $a_{ij}=1$ if and only if $v_i$ is adjacent
to $v_j$. Then $A$ is symmetric and $(1,\ldots,1) \in\mathbb{R}^d$
is an eigenvector of $A$ with eigenvalue $k$. The set of
eigenvalues of $A$ is $\{k,r_{\!{}_+},r_{\!{}_-}\}$ with
$r_{\!{}_-}< 0\leq r_{\!{}_+}\leq k$ and
\begin{equation}
\label{roots}
 r_{\!{}_\pm} =\frac{1}{2}\left[\,\lambda -\mu \pm \sqrt{(\lambda
-\mu)^2 +4(k-\mu)}\,\,\right],
\end{equation} implying that
\begin{equation}
\label{quadeq1} (A-r_{\!{}_+}I_d)(A-r_{\!{}_-}I_d)
=\frac{(k-r_{\!{}_+})(k-r_{\!{}_-})}{d}J_d\,,
\end{equation} where $J_d$ is the $d\times d$ matrix whose entries
are equal to 1. Given parameters $(d,k\lambda,\mu)$ such that
$(d,k,\lambda,\mu)\neq (4\mu +1,2\mu,\mu -1,\mu)$, one can show
that $r_{\!{}_\pm} \in \mathbb{Z}$ (cf.\ \cite{LW}, Theorem 21.1).
In fact, if $(d,k,\lambda,\mu)=(4\mu +1,2\mu,\mu -1,\mu)$, then
non-integral values of $r_{\!{}_\pm}$ can occur; for instance, the
{\it Paley} graph $P(5)$ (see \cite{LW}, Example 21.3) has
parameters $(5,2,0,1)$ and eigenvalues
$r_{\!{}_\pm}=-\frac{1}{2}\pm\frac{1}{2}\sqrt{5}$. However, the
Paley graph $P(4\mu +1)\in {{\rm SRG}(4\mu +1,2\mu,\mu -1,\mu)}$,
where $4\mu +1=p^{2n}$ with $p$ an odd prime and $n\in
\mathbb{N}$, has
integer eigenvalues $r_{\!{}_\pm}=\tfrac{1}{2}(1\pm p^n)$.\\

Many strongly regular graphs stem from geometry. The most
classical example is offered by the configuration of 27 lines on a
cubic surface.

\begin{example}
\label{schlaefligraph} Given a non-singular cubic surface
$X\subset \mathbb{P}^3$, let $\mathcal{L}$ be the intersection
graph of the 27 lines that are contained in $X$. The configuration
of the 27 lines is completely described by the 36 {\it Sch\"{a}fli
double-sixes}, i.e., pairs
$M:=\bigl(\{a_1,\ldots,a_6\},\{b_1,\ldots,b_6\}\bigr)$ of sets of
6 skew lines on $X$ such that each line from one set is skew to a
unique line from the other set. Fix a double-six $M$; in matrix
notation we may write
$$
M=\left( \begin{matrix} a_1 & a_2 & a_3 & a_4 & a_5 & a_6 \\
b_1 & b_2 & b_3 & b_4 & b_5 & b_6 \end{matrix} \right)
$$ such that two lines meet if and only they are in different rows
and columns. The remaining 15 lines on $X$ are the
$c_{ij}:=a_ib_j\cap a_jb_i$ with $i\neq j$, where $a_ib_j$ is the
plane spanned by the lines $a_i$ and $b_j$. In this notation the
36 double-sixes are $M$, the 15 $M_{i,j}$'s and the 20
$M_{i,j,k}$'s, where the double-sixes $M_{i,j}$ and $M_{i,j,k}$
are respectively given by
$$
\left( \begin{matrix} a_i & b_i & c_{jk} & c_{jl} & c_{jm} & c_{jn} \\
a_j & b_j & c_{ik} & c_{il} & c_{im} & c_{in} \end{matrix}
\right),\hskip0.25cm
\left( \begin{matrix} a_i & a_j & a_k & c_{mn} & c_{ln} & c_{ln} \\
c_{jk} & c_{ik} & c_{ij} & b_l & b_m & b_n \end{matrix} \right).
$$ It follows that the {\it Schl\"{a}fli} graph $\mathcal{L}$ is
in ${\rm SRG}(27,10,1,5)$ (the unique member of this set). It is
easily seen that the stabilizer of a double-six is a subgroup of
${\rm Aut}(\mathcal{L})$ of index 36, isomorphic to $S_6\times
\mathbb{Z}/2\mathbb{Z}$; consequently $\#({\rm
Aut}(\mathcal{L}))=6\,!\cdot 2\cdot 36$. In fact, ${\rm
Aut}(\mathcal{L})$ is isomorphic to the Weyl group $W(E_6)$:
Consider the Dynkin diagram for $E_6$

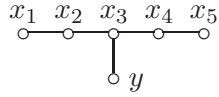
\begin{center}
\setlength{\unitlength}{3mm}
\begin{picture}(30,5)(-10,-3)

\multiput(0,0)(2,0){5}{$\circ$}
\multiput(0.5,.4)(2,0){4}{\line(1,0){1.6}} \put(4,-2){$\circ$}
\put(4.3,-1.4){\line(0,1){1.5}}

\put(-0.3,1){$x_1$}\put(1.7,1){$x_2$}\put(3.7,1){$x_3$}
\put(5.7,1){$x_4$}\put(7.7,1){$x_5$}\put(5,-2){$y$}

\end{picture}
\end{center}

\noindent By definition $W(E_6)$ is generated by the reflections
$s_1,\ldots,s_5,s$ associated to the simple roots
$x_1,\ldots,x_5,y$. One shows that there is a surjective
homomorphism $W(E_6)\rightarrow {\rm Aut}(\mathcal{L})$ sending
$s_i$ (resp.\ $s$) to the transformation that interchanges the
rows of the double-six $M_{i,i+1}$ (resp.\ $M_{1,2,3}$) (cf.\
\cite{Ma}, sections 25,26). Then recall that $W(E_6)$ is of order
$6\,!\cdot 2\cdot 36$. Under this apparent isomorphism the 27
lines on $X$ correspond to the 27 fundamental weights of $E_6$.
Since $W(E_6)$ acts transitively on these weights, we may consider
${\rm Aut}(\mathcal{L})$ as a transitive subgroup of $S_{27}$.
\end{example}

\noindent \textbf{Additional properties} (of strongly regular
graphs):
\begin{itemize}
\item[1)] if $\mathcal{G}$ is disconnected, then $\mathcal{G}$ is
the disjoint union of $m>1$ copies of the {\it complete} graph
$K_{k+1}$ with adjacency matrix $J_{k+1} -I_{k+1}$, where
$I_{k+1}$ is the identity matrix. So
$(d,k,\lambda,\mu)={\bigl(m(k+1),k,k-1,0\bigr)}$ and $A$ has
exactly two distinct eigenvalues: $k$ $(=r_{\!{}_+})$ with
multiplicity $m$ and $r_{\!{}_-}=-1$ with multiplicity $d-m$;
\item[2)] if $\mathcal{G}$ is connected and $\mathcal{G} \neq
K_{k+1}$, then $k\neq r_{\!{}_\pm}$ and $k$ is a simple
eigenvalue;
\item[3)] $\mathcal{G}$ and its complement
$\overline{\mathcal{G}}\in {\rm
SRG}(d,d-k-1,d-2k+\mu-2,d-2k+\lambda)$ are connected if and only
if $0<\mu<k<d-1$, in which case $\mathcal{G}$ is said to be {\it
non-trivial}.
\end{itemize} Let ${\rm Aut}(A)$ be the stabilizer of $A$
under the natural operation of $S_d$ on $d\times d$ integer
matrices by
$(a_{ij})\overset\sigma\mapsto(a_{\sigma(i)\sigma(j)})$ and
observe that ${\rm Aut}(A)$ coincides with ${\rm
Aut}(\mathcal{G})$. For disconnected $\mathcal{G}$ it is easily
seen that ${\rm Aut}(A)$ is a transitive subgroup of $S_d$,
implying that disconnected strongly regular graphs are {\it
vertex-transitive}. In practice it turns out that
vertex-transitivity is quite rare among non-trivial strongly
regular graphs, although most of the sets ${\rm
SRG}(d,k,\lambda,\mu)$ of non-trivial strongly regular graphs have
a vertex-transitive member.

\begin{definition}
\label{prymmatrix} Let $A\in \mathbb{Z}^{d\times d}$ be a
symmetric matrix with transitive stabilizer group ${\rm Aut}(A)$.
Then $A$ is a {\it Prym} matrix if there exist integers
$k,r_{\!{}_+},r_{\!{}_-}$ with $r_{\!{}_+}>r_{\!{}_-}$ such that
equation \eqref{quadeq1} holds. Further, if $A$ is a Prym matrix
and $m\in \mathbb{N}$, then $A^{\oplus m}:=\bigoplus_{i=1}^{m}A$
is a {\it repeated} Prym matrix.
\end{definition}

\noindent \textbf{Remarks.} Suppose that $A\in \mathbb{Z}^{d\times
d}$ is a symmetric matrix for which there exist integers
$k,r_{\!{}_+},r_{\!{}_-}$ with $r_{\!{}_+}>r_{\!{}_-}$ such that
equation \eqref{quadeq1} holds. Decomposing $\mathbb{R}^d$ into
eigenspaces of $A$ we may assume that $A$ takes diagonal form.
Then $J_d$ simultaneously transforms into the diagonal matrix
${\rm diag}(d,0,\ldots,0)$, implying that $(1,\ldots,1)\in
\mathbb{R}^d$ is an eigenvector of $A$ with, say, eigenvalue
$\eta$. Hence $A$ has eigenvalues $\eta,r_{\!{}_+},r_{\!{}_-}$ and
by \eqref{quadeq1} we have
$(\eta-r_{\!{}_+})(\eta-r_{\!{}_-})=(k-r_{\!{}_+})(k-r_{\!{}_-})$,
that is, $\eta=k$ or $\eta=r_{\!{}_+} +r_{\!{}_-} -k$. Further, if
$\eta\neq r_{\!{}_\pm}$, then $\eta$ is simple. A $d\times d$ Prym
matrix $A$ therefore has integer eigenvalues
$k,r_{\!{}_+},r_{\!{}_-}$ with $r_{\!{}_+}>r_{\!{}_-}$ such that
\eqref{quadeq1} holds and $k$ belongs to the eigenvector
$(1,\ldots,1)$ of $A$. Clearly, if $A$ is such a matrix, then for
any $m\in \mathbb{N}$ the repeated Prym matrix $A^{\oplus m}$ has
the same eigenvalues $k,r_{\!{}_+},r_{\!{}_-}$ and satisfies the
equation
\begin{equation}
\label{quadeq2} (A^{\oplus m}-r_{\!{}_+}I_{md})(A^{\oplus
m}-r_{\!{}_-}I_{md})=\frac{(k-r_{\!{}_+})(k-r_{\!{}_-})}{d}J_d^{\oplus
m}.
\end{equation} Moreover, it is immediately seen that the repeated
Prym matrix $A^{\oplus m}$ has transitive stabilizer ${\rm Aut}
(A^{\oplus m})$. Throughout the paper the eigenvalues of a Prym
matrix $A$ will be denoted by $k,r_{\!{}_+},r_{\!{}_-}$, where $k$
belongs to the eigenvector $(1,\ldots,1)$ of $A$ and
$r_{\!{}_+}>r_{\!{}_-}$.\\

Many of the known constructions for Prym varieties rely on the
definition of a reduced symmetric divisor correspondence for a
curve. These correspondences can often be related to Prym matrices
whose entries are in the set $\{0,1\}$. Hence it is useful to have
a characterization for such matrices:

\begin{proposition}
\label{srg} Assume that $A=(a_{ij})\in \{0,1\}^{d\times d}$ is a
Prym matrix with $a_{11}=0$. Then $A$ is the adjacency matrix of a
graph $\mathcal{G}\in {\rm SRG}(d,k,\lambda,\mu)$ with
$\lambda=k+r_{\!{}_+}r_{\!{}_-}+r_{\!{}_+}+r_{\!{}_-}$ and
$\mu=k+r_{\!{}_+}r_{\!{}_-}$.
\end{proposition}

\proof By transitivity we have $a_{ii}=0$ for all $i=1,\ldots,d$,
hence $A$ is the adjacency matrix of a regular graph $\mathcal{G}$
of degree $k$ with $d$ vertices. If indeed we have $\mathcal{G}\in
{\rm SRG}(d,k,\lambda,\mu)$ for some integers $\lambda$ and $\mu$,
then solving equation \eqref{roots} for $\lambda,\mu$ yields
$\lambda=k+r_{\!{}_+}r_{\!{}_-}+r_{\!{}_+}+r_{\!{}_-}$ and
$\mu=k+r_{\!{}_+}+r_{\!{}_-}$. To show that $\mathcal{G}$ is
strongly regular we may assume by transitivity that $\mathcal{G}$
is the disjoint union of finitely many copies of a connected graph
$\mathcal{G}_c$ with adjacency matrix~$A_c$. If $k\neq
r_{\!{}_\pm}$, then $k$ is simple, so $A=A_c$ and it follows that
$\mathcal{G}$ is strongly regular (cf.\ \cite{LW}, the remark at
the bottom of p.\ 265). Hence assume that $A$ has just two
distinct eigenvalues, i.e., $k=r_{\!{}_+}$. The complementary
graph $\overline{\mathcal{G}}$ has adjacency matrix $J_d-I_d-A$
which is Prym and has eigenvalues $k'=d-k-1$,
$r'_{\!{}_+}=-r_{\!{}_-}-1$ and $r'_{\!{}_-}=-k-1$. So if
$r_{\!{}_-}\neq k-d$, then $k'\neq r'_{\!{}_\pm}$, that is,
$\overline{\mathcal{G}}$ is a non-complete connected strongly
regular graph and so $\mathcal{G}_c=K_{k+1}$. Finally, suppose
that $r_{\!{}_-}=k-d$. By equation \eqref{quadeq1} we have
$A^2=(2k-d)A-k(k-d)I_d$, hence if $v_h,v_i,v_j$ are three distinct
vertices of $\mathcal{G}$ such that $v_h$ is adjacent to $v_i$ and
$v_i$ is adjacent to $v_j$, then the relation
$1=a_{hi}\leqq\sum_{l=1}^d a_{hl}a_{lj}=(2k-d)a_{hj}$ implies that
$v_h$ is adjacent to $v_j$. Consequently
$\mathcal{G}=K_{k+1}$.\qed

\begin{example}
\label{latticegraphs} Among the set of strongly regular graphs
there are some distinguished families of non-trivial vertex
transitive graphs. One such family is that of {\it lattice}
graphs; for $n\geq 3$, the lattice graph $L_2(n)$ is the graph
with vertex set $\{1,\ldots,n\}^2$ such that two distinct vertices
$(i,j)$ and $(l,m)$ are adjacent if and only if $i=l$ or $j=m$.
Clearly, $S_n\times S_n$ is a transitive subgroup of ${\rm
Aut}\bigl(L_2(n)\bigr)$, hence the adjacency matrix of $L_2(n)$
and that of its complement $\overline{L_2(n)}$ are Prym. Another
subgroup of ${\rm Aut}\bigl(L_2(n)\bigr)$ is $S_2$; it permutes
the coordinates of the vertices of $L_2(n)$. In fact, it is
well-known that the automorphism group of $L_2(n)\in {\rm
SRG}\bigl(n^2,2(n-1),n-2,2\bigr)$ is equal to the semi-direct
product $S_2\rtimes (S_n\times S_n)$. With reference to future
examples (e.g., Examples \ref{lattice} and \ref{twisted}) we shall
characterize those $\varphi\in {\rm Aut}\bigl(L_2(n)\bigr)$ for
which each vertex $(i,j)$ of $\overline{L_2(n)}$ (resp.\ $L_2(n)$
for $n$ odd) is non-adjacent to $\varphi(i,j)$. It is easy to
check that:
\begin{itemize}
\item[i)] Each vertex $(i,j)$ of $\overline{L_2(n)}$ is
non-adjacent to $\varphi(i,j)$ if and only if
$\varphi=(\sigma,\tau)$ with $\sigma,\tau\in S_n$ and $\sigma=(1)$
or $\tau=(1)$.
\item[ii)] Assume that $n$ is odd and $\varphi$ is
a reflection. Then each vertex $(i,j)$ of $L_2(n)$ is non-adjacent
to $\varphi(i,j)$ if and only if
$\varphi=(\sigma,\sigma^{-1})\circ t$ with $\sigma\in S_n$ and
$t=(1\hskip0.25cm2)\in S_2$.
\end{itemize} Note moreover that $L_2(3)$ and
$\overline{L_2(3)}$ are isomorphic; if we identify the set
$\{1,\ldots,n\}$ with the group $\mathbb{Z}/3\mathbb{Z}$, then the
matrix $\bigl(
\begin{smallmatrix} 1 & -1 \\ 1 & 1 \end{smallmatrix} \bigr)$
defines an isomorphism of graphs $L_2(3) \overset \sim \rightarrow
\overline{L_2(3)}$.
\end{example}

\begin{example}
\label{latinsquaregraphs} Another non-trivial family is that of
{\it Latin square} graphs; for $n\geq 3$, the Latin square graph
$L_3(n)\in {\rm SRG}(n^2,3(n-1),n,6)$ is the graph with vertex set
$(\mathbb{Z}/n\mathbb{Z})^2$ such that two distinct vertices
$(i,j)$ and $(l,m)$ are adjacent if and only if $i=l$ or $j=m$ or
$i+j=l+m$. We identify three subgroups of ${\rm
Aut}\bigl(L_3(n)\bigr)$; to begin with, the diagonal action of
${\rm Aut}(\mathbb{Z}/n\mathbb{Z})$ on
$(\mathbb{Z}/n\mathbb{Z})^2$ induces an embedding ${\rm
Aut}(\mathbb{Z}/n\mathbb{Z})\hookrightarrow {\rm
Aut}\bigl(L_3(n)\bigr)$. Similarly, $(\mathbb{Z}/n\mathbb{Z})^2$
induces a subgroup of ${\rm Aut}\bigl(L_3(n)\bigr)$ by
translation, thus implying that $L_3(n)$ is vertex-transitive.
Finally, the subgroup $S:=\langle s,t\rangle\subset{\rm
Aut}\bigl((\mathbb{Z}/n\mathbb{Z})^2\bigr)$ with $s$ and $t$
respectively given by the matrices $\bigl(
\begin{smallmatrix} -1 & -1 \\ 0 & 1 \end{smallmatrix} \bigr)$ and
$\bigl(
\begin{smallmatrix} 0 & 1 \\ 1 & 0 \end{smallmatrix} \bigr)$, is
immediately seen to be a subgroup of ${\rm
Aut}\bigl(L_3(n)\bigr)$. Note that sending $s\mapsto
(1\hskip0.25cm 2)$ and $t\mapsto (1\hskip0.25cm 3)$ induces an
isomorphism $S\overset \sim \rightarrow S_3$. Clearly, the actions
of ${\rm Aut}(\mathbb{Z}_n)$ and $S$ commute. It can be shown that
${\rm Aut}\bigl(L_3(n)\bigr)$ coincides with the semi-direct
product $(\mathbb{Z}/n\mathbb{Z})^2\rtimes \bigl(S\times {\rm
Aut}(\mathbb{Z}/n\mathbb{Z})\bigr)$.
\end{example}

\end{section}

\begin{section}{Prym varieties of a Galois covering}
\label{galoisconstruction}

Given a Prym matrix $A$, we describe a method that associates
certain abelian varieties to a finite Galois coverings of
$\mathbb{P}^1$ whose Galois group is represented as a transitive
subgroup of ${\rm Aut}(A^{\oplus m})$, $m\geq 1$. We assume from
now on that $A$ is a $d\times d$ Prym matrix with eigenvalues
$k,r_{\!{}_+},r_{\!{}_-}$.

\begin{definition}
\label{prymdata} Consider $\mathcal{P}=(A^{\oplus m},r,\pi,\phi)$
for $m\geq 1$, where $r\in\{r_{\!{}_+},r_{\!{}_-}\}$, $\pi$ is a
finite Galois covering of $\mathbb{P}^1$ and $\phi:{\rm
Gal}(\pi)\rightarrow S_{md}$ is a transitive representation. Then
$\mathcal{P}$ is said to represent {\it Prym data} if $\phi(G)$ is
a subgroup of ${\rm Aut}(A^{\oplus m})$.
\end{definition}

For instance, assume that we have a finite subset
$B=\{b_1,\ldots,b_n\}$ of $\mathbb{P}^1$ along with non-trivial
permutations $\sigma_1,\ldots,\sigma_n\in {\rm Aut}(A^{\oplus m})$
such that $\sigma_1\cdot\ldots\cdot\sigma_n=(1)$ and $G:=\langle
\sigma_1,\ldots,\sigma_n\rangle$ is a transitive subgroup of
$S_{md}$. Let $\Sigma_i$ be the conjugation class of $\sigma_i$ in
$G$. According to Riemann's Existence Theorem (RET), the number of
equivalence classes of Galois coverings of $\mathbb{P}^1$ of
ramification type $\mathcal{R}:=[G,B,\{\Sigma_i\}_{i=1}^n]$ is
equal to the number of sets $\{(g\tau_1 g^{-1},\ldots,g\tau_n
g^{-1})|\,g\in G\}$ with $\tau_i\in \Sigma_i$ such that
$\tau_1\cdot\ldots\cdot\tau_n$ is trivial (cf.\ \cite{Vo}, p.\
37). Hence, let $\pi$ be a Galois covering of type $\mathcal{R}$.
Clearly, if $\phi:{\rm Gal}(\pi)\overset\sim\rightarrow G$ is a
group isomorphism, then $(A^{\oplus m},r,\pi,\phi)$ represents
Prym data. In this way, varying the continuous parameters
$b_1,\ldots,b_n$, we obtain finitely many smooth $n$-dimensional
families of Galois coverings with associated Prym data.\\

\noindent \textbf{The construction.} Let $\mathcal{P}=(A^{\oplus
m},r,\pi,\phi)$ be Prym data for a given Galois covering
$\pi:X\rightarrow \mathbb{P}^1$. We are going to define two
symmetric divisor correspondences, one for $X$ and one for the
quotient curve $C:=X/H$, where $H\subset G:={\rm Gal}(\pi)$ is the
stabilizer of the letter 1 with respect to $\phi$. Let
$(\,\,,\,)_r:\mathbb{R}^{md}\times\mathbb{R}^{md}\rightarrow
\mathbb{Q}$ be the symmetric bilinear form canonically associated
to the matrix $A^{\oplus m}-rI_{md}$. Given the standard basis
$e_1,\ldots,e_{md}$ of $\mathbb{R}^{md}$, consider the permutation
representation of $G$ on $\mathbb{R}^{md}$ induced by
$g:e_i\mapsto e_{g(i)}$. For $g\in G$, denote $\hat{g}=HgH$. Then
$(\,\,,\,)_r$ is immediately seen to be $G$-invariant and
$(g_1e_1,e_1)_r=(g_2e_1,e_1)_r$ for all $g_1,g_2\in G$ such that
$\hat{g}_1=\hat{g}_2$. Let $\alpha:X\rightarrow C$ and
$p:C\rightarrow \mathbb{P}^1$ be the canonical mappings. For each
$g\in G$ take the graph $\Gamma_g=({\rm id}_X,g)(X)$ of $g$ and
let $\hat{\Gamma}_g=(\alpha\times \alpha)(\Gamma_g)$ be its
reduced image in $C\times C$. Assume that $B$ is the branch locus
of $\pi$ and put $C_0=p^{-1}(\mathbb{P}^1\setminus B)$. Observe
that
$$\hat{\Gamma}_{g_1}\cap\hat{\Gamma}_{g_2}\cap
(C_0\times C_0)\neq\emptyset\Longleftrightarrow\hat{g}_1=\hat{g}_2
\Longleftrightarrow\hat{\Gamma}_{g_1}=\hat{\Gamma}_{g_2}\,,
$$ for all $g_1,g_2\in G$. Hence we have divisor correspondences
$D_{\!{}_\mathcal{P}}$ on $X\times X$ and
$\hat{D}_{\!{}_\mathcal{P}}$ on $C\times C$, given by
$$
D_{\!{}_\mathcal{P}}=\sum_{g\in G}(ge_1,e_1)_r
\Gamma_g\,,\hskip0.5cm \hat{D}_{\!{}_\mathcal{P}}=\sum_{\hat{g}\in
\hat{G}}(ge_1,e_1)_r\hat{\Gamma}_g
$$ with $\hat{G}=H\backslash G/H$. Note that $D_{\!{}_\mathcal{P}}$
and $\hat{D}_{\!{}_\mathcal{P}}$ are symmetric as
$(ge_1,e_1)_r=(g^{-1}e_1,e_1)_r$ for all $g\in G$.

\begin{definition}
\label{prym1} Let $\gamma_{\!{}_\mathcal{P}}$ on $J(X)$ (resp.\
$\hat{\gamma}_{\!{}_\mathcal{P}}$ on $J(C)$) be the endomorphism
canonically associated to $D_{\!{}_\mathcal{P}}$ (resp.\
$\hat{D}_{\!{}_\mathcal{P}}$). Then we call
$Z_{\!{}_\mathcal{P}}={\rm im}\,\gamma_{\!{}_\mathcal{P}}$ (resp.\
$\hat{Z}_{\!{}_\mathcal{P}}={\rm
im}\,\hat{\gamma}_{\!{}_\mathcal{P}}$) the Prym variety of $X$
(resp.\ $C$) associated to $\mathcal{P}$.
\end{definition}

\noindent \textbf{Remark.} Prym data can be seen as an
`ornamented' covering (i.e., a covering with additional data),
where the ornamentation is such that we can define divisorial
correspondences and Prym varieties.\\

For $q\in \mathbb{P}^1\setminus B$ an identification
$\pi^{-1}(q)\leftrightarrow G$ is called a {\it Galois labelling}
of the fiber $\pi^{-1}(q)$ if the action of $G$ on the fiber is
compatible with its action on itself via multiplication on the
left. Moreover, a Galois labelling of $\pi^{-1}(q)$ induces a
Galois labelling $p^{-1}(q)\leftrightarrow H\backslash G$ of the
fiber of $p:C\rightarrow \mathbb{P}^1$ over $q$. Denoting
$\overline{g}=Hg$ for $g\in G$, we have:

\begin{lemma}
\label{identity} Given $q\in \mathbb{P}^1\setminus B$, take a
Galois labelling $p^{-1}(q)\leftrightarrow H\backslash G$. Then,
for any $\sigma\in G$, the restriction of
$\hat{D}_{\!{}_\mathcal{P}}$ to $\{\overline{\sigma}\}\times C$ is
given by the identity
$$
\hat{D}_{\!{}_\mathcal{P}}(\overline{\sigma})=\sum_{\overline{g}\in
H\backslash G}(ge_1,e_1)_r\,\overline{g\sigma}\,.
$$
\end{lemma}

\proof Because $(\overline{\sigma},\overline{g\sigma})\in
\hat{\Gamma}_{\!f}\Longleftrightarrow \hat{f}=\hat{g}$, for all
$f,g\in G$, the point $(\overline{\sigma},\overline{g\sigma})$
appears in $\hat{D}_\mathcal{P}$ with multiplicity
$(ge_1,e_1)_r$.\qed\\

Since $H$ is the stabilizer of a point, the transitivity of the
representation $\phi:G\rightarrow S_{md}$ implies that there
exists a bijection $H\backslash G\rightarrow\{1,\ldots,md\}$.
Hence the Galois labelling $p^{-1}(q)\leftrightarrow H\backslash
G$ induces a labelling $\{y_1,\ldots,y_{md}\}$ of $p^{-1}(q)$ such
that $\overline{g}$ corresponds to $y_{g^{-1}(1)}$. The formula in
the preceding lemma now turns into
\begin{equation}
\label{matrixidentity}
\hat{D}_{\!{}_\mathcal{P}}(y_i)=\sum_{j=1}^{md}(e_i,e_j)_r\,
y_j=-ry_i+\sum_{j=1}^{md}({}^t\!e_jA^{\oplus m}e_i)\,y_j\,,
\end{equation} for all $i=1,\ldots,md$.

\begin{proposition}
\label{isogeny} Let $\mathcal{P}=(A^{\oplus m},r,\pi,\phi)$ be
Prym data with $\pi:X\rightarrow \mathbb{P}^1$. Then the Prym
varieties $Z_{\!{}_\mathcal{P}}$ of $X$ and
$\hat{Z}_{\!{}_\mathcal{P}}$ of $C:=X/H$ are isogenous. Using the
notation $\mathcal{P}'=(J_d^{\oplus m},0,\pi,\phi)$ and
$s=r_{\!{}_+}+r_{\!{}_-}$, we have the following quadratic
equations for the endomorphisms $\hat{\gamma}_{\!{}_\mathcal{P}}$
on $J(C)$ and $\gamma_{\!{}_\mathcal{P}}$ on $J(X)$:
$$
\hat{\gamma}_{\!{}_\mathcal{P}}
\bigl(\hat{\gamma}_{\!{}_\mathcal{P}}+(2r-s){\rm id}_{J(C)}\bigr)=
\frac{(k-r_{\!{}_+})(k-r_{\!{}_-})}{d}\,
\hat{\gamma}_{\!{}_{\mathcal{P}'}}
$$ and
$$\gamma_{\!{}_\mathcal{P}}
\bigl(\gamma_{\!{}_\mathcal{P}}+\#(H) (2r-s){\rm id}_{J(X)}\bigr)=
\#(H)\frac{(k-r_{\!{}_+})(k-r_{\!{}_-})}{d}\,
\gamma_{\!{}_{\mathcal{P}'}}.
$$
\end{proposition}

\proof Write $\alpha^*:J(C)\rightarrow J(X)$ for the map induced
by $\alpha$ and write $N_{\alpha}:J(X)\rightarrow J(C)$ for the
norm map. Abusing notation, we define homomorphisms
$\alpha^*:\mathbb{Z}[H\backslash G]\rightarrow \mathbb{Z}[G]$ and
$N_{\alpha}:\mathbb{Z}[G]\rightarrow \mathbb{Z}[H\backslash G]$,
respectively induced by $\overline{g}\mapsto \sum_{h\in H}hg$ and
$g\mapsto \overline{g}$. As a direct consequence of Lemma
\ref{identity} we have $N_{\alpha}D_{\!{}_\mathcal{P}}=
\#(H)\hat{D}_{\!{}_\mathcal{P}}N_{\alpha}$. Hence
$N_{\alpha}\gamma_{\!{}_\mathcal{P}}=
\#(H)\hat{\gamma}_{\!{}_\mathcal{P}}N_{\alpha}$ and
$\gamma_{\!{}_\mathcal{P}}=\alpha^*
\hat{\gamma}_{\!{}_\mathcal{P}}N_{\alpha}$. The first identity
implies that the restricted mapping
$N_{\alpha}:Z_{\!{}_\mathcal{P}}\rightarrow
\hat{Z}_{\!{}_\mathcal{P}}$ is surjective, while the two
identities combined imply that ${\rm
dim\,im}\,Z_{\!{}_\mathcal{P}}={\rm
dim\,im}\,\hat{Z}_{\!{}_\mathcal{P}}$. Therefore the restriction
of $N_{\alpha}$ to $Z_{\!{}_\mathcal{P}}$ is an isogeny.

Applying equation \eqref{matrixidentity} to $\mathcal{P}'$ we
obtain $\hat{D}_{\!{}_{\mathcal{P}'}}(y_i)=
\sum_{j=1}^{md}({}^t\!e_jJ_d^{\oplus m}e_i)\,y_j$, for all
$i=1,\ldots,md$. Hence equations \eqref{quadeq2} and
\eqref{matrixidentity} imply
$$
\hat{D}_{\!{}_\mathcal{P}}
\bigl(\hat{D}_{\!{}_\mathcal{P}}(y)\bigr)+
(2r-s)\hat{D}_{\!{}_\mathcal{P}}(y)=
\frac{(k-r_{\!{}_+})(k-r_{\!{}_-})}{d}
\,\hat{D}_{\!{}_{\mathcal{P}'}}(y)
$$ for all $y$ in a fiber of
$p:C\rightarrow \mathbb{P}^1$ over a point outside the branch
locus of~$\pi$. Let $\hat{\delta}$ (resp.\ $\delta$) denote the
difference between the expressions on the left and the righthand
side of the first (resp.\ second) identity of the proposition. We
already know that $\hat{\delta}=0$. Hence, as
$\gamma_{\!{}_{\mathcal{P}'}}= \alpha^*
\hat{\gamma}_{\!{}_{\mathcal{P}'}}N_{\alpha}$, we have $\delta=
\#(H)\alpha^*\hat{\delta}N_{\alpha}=0$.\qed  

\begin{proposition}
\label{complement} In the special case $m=1$ the endomorphisms
$\hat{\gamma}_{\!{}_{\mathcal{P}'}}$ and
$\gamma_{\!{}_{\mathcal{P}'}}$ vanish and, using the notation
$\mathcal{P}^{{}_\pm}=(A,r_{\!{}_\pm},\pi,\phi)$, we find that
$\hat{Z}_{\!{}_{\mathcal{P}^+}}$ and
$\hat{Z}_{\!{}_{\mathcal{P}^-}}$ are complementary subvarieties of
$J(C)$. Moreover, $\hat{Z}_{\!{}_{\mathcal{P}^\pm}}={\rm
ker}(\hat{\gamma}_{\!{}_{\mathcal{P}^\mp}})_{\!{}_0}$ and
$Z_{\!{}_{\mathcal{P}^\pm}}={\rm ker}
\bigl(\gamma_{\!{}_{\mathcal{P}^\pm}}\pm\#(H)(r_{\!{}_+}-r_{\!{}_-}){\rm
id}_{J(X)}\bigr)_{\!{}_0}$.
\end{proposition}

\proof Let $q\in \mathbb{P}^1\setminus B$. By definition of
$D_{\!{}_{\mathcal{P}'}}$ we have
$D_{\!{}_{\mathcal{P}'}}(x)=\pi^*\pi(x)$ for all $x\in
\pi^{-1}(q)$, while Lemma \ref{identity} implies that
$\hat{D}_{\!{}_{\mathcal{P}'}}(y)=p^*p(y)$ for all $y\in
p^{-1}(q)$. Hence $\gamma_{\!{}_{\mathcal{P}'}}$ and
$\hat{\gamma}_{\!{}_{\mathcal{P}'}}$ vanish. As
$\hat{\gamma}_{\!{}_{\mathcal{P}^-}}=
\hat{\gamma}_{\!{}_{\mathcal{P}^+}}+(r_{\!{}_+}-r_{\!{}_-}){\rm
id}_{J(C)}$, we note that $\hat{Z}_{\!{}_{\mathcal{P}^-}}={\rm
im}\bigl(\hat{\gamma}_{\!{}_{\mathcal{P}^+}}+(r_{\!{}_+}-r_{\!{}_-}){\rm
id}_{J(C)}\bigr)$. Using standard arguments one shows that
$\varepsilon_{\!{}_+}:=\frac{-1}{r_{\!{}_+}-r_{\!{}_-}}
\hat{\gamma}_{\!{}_{\mathcal{P}^+}}$ and
$\varepsilon_{\!{}_-}:=\frac{1}{r_{\!{}_+}-r_{\!{}_-}}
\bigl(\hat{\gamma}_{\!{}_{\mathcal{P}^+}}+(r_{\!{}_+}-r_{\!{}_-}){\rm
id}_{J(C)}\bigr)$ are symmetric idempotents in ${\rm
End}_\mathbb{Q}\bigl(J(C)\bigr)$. Since $\varepsilon_{\!{}_+}={\rm
id}_{J(C)}-\varepsilon_{\!{}_-}$, it follows that
$Z_{\!{}_{\mathcal{P}^+}}$ and $Z_{\!{}_{\mathcal{P}_-}}$ are
complementary subvarieties of $J(C)$ and ${\rm
dim}\,Z_{\!{}_{\mathcal{P}^+}}+{\rm
dim}\,Z_{\!{}_{\mathcal{P}^-}}=g(C)$. As a consequence of
Proposition \ref{isogeny} we have
$\hat{Z}_{\!{}_{\mathcal{P}^\pm}}\subset {\rm
ker}(\hat{\gamma}_{\!{}_{\mathcal{P}^\mp}})_{\!{}_0}$. Moreover,
by \cite{BL}, Proposition 5.1.1 the analytic representation
$\varrho_a(\hat{\gamma}_{\!{}_{\mathcal{P}^+}})\in {\rm
End}\bigl(H^0(C,\omega_{C})\bigr)$ of
$\hat{\gamma}_{\!{}_{\mathcal{P}^+}}$ is self-adjoint with
eigenvalues $r_{\!{}_-}-r_{\!{}_+}$ and 0 with respect to the
Riemann form $c_1(\Theta_{C})$, where $\Theta_C$ is the canonical
polarization of $J(C)$. Since ${\rm
dim\,ker}(\hat{\gamma}_{\!{}_{\mathcal{P}^\pm}})_{\!{}_0}= {\rm
dim\,ker}\,\varrho_a(\hat{\gamma}_{\!{}_{\mathcal{P}^\pm}})$, we
have ${\rm
dim\,ker}(\hat{\gamma}_{\!{}_{\mathcal{P}^+}})_{\!{}_0}+{\rm
dim\,ker}(\hat{\gamma}_{\!{}_{\mathcal{P}^-}})_{\!{}_0}=g(C)$ and
therefore $\hat{Z}_{\!{}_{\mathcal{P}^\pm}}={\rm
ker}(\hat{\gamma}_{\!{}_{\mathcal{P}^\mp}})_{\!{}_0}$. Similarly,
the remaining assertion follows from the fact that
$Z_{\!{}_{\mathcal{P}^\pm}}$ and ${\rm im}
\bigl(\gamma_{\!{}_{\mathcal{P}^\pm}}
\pm\#(H)(r_{\!{}_+}-r_{\!{}_-}){\rm id}_{J(X)}\bigr)$ are
complementary subvarieties of $J(X)$. \qed\\

Finally, we note that the Prym data $\mathcal{P}=(A,r,\pi,\phi)$
associated to the adjacency matrix $A$ of a strongly regular graph
$\mathcal{G}$ and the Prym data
$\overline{\mathcal{P}}=(J_d-I_d-A,-r-1,\pi,\phi)$ associated to
the adjacency matrix $J_d-I_d-A$ of the complementary graph
$\overline{\mathcal{G}}$ yield the same Prym varieties because of
the identities
$\gamma_{\!{}_\mathcal{P}}+\gamma_{\!{}_{\overline{\mathcal{P}}}}=0$
and $\hat{\gamma}_{\!{}_\mathcal{P}}
+\hat{\gamma}_{\!{}_{\overline{\mathcal{P}}}}=0$.

\end{section}

\begin{section}{Prym varieties of a non-Galois covering}
\label{nongalois}

We now shift our attention to non-Galois coverings of
$\mathbb{P}^1$. Given a (repeated) Prym matrix $A^{\oplus m}$ and
a special covering $p:C\rightarrow \mathbb{P}^1$ of degree~$md$
along with an appropriate {\it labelling class} (to be defined
below), we shall define a symmetric divisor correspondence on
$C\times C$ which then serves to obtain a pair of Prym varieties
in $J(C)$.\\

To define labelling classes, let $p:C\rightarrow \mathbb{P}^1$ be
a covering of degree $n$ with branch locus $B$. Given a point
$q\in \mathbb{P}^1\setminus B$, a labelling $\{x_1,\ldots,x_n\}$
of the fiber $p^{-1}(q)$ induces a bijection
$\nu:p^{-1}(q)\rightarrow\{1,\ldots,n\}$ sending $x_i$ to $i$. For
$q_j\in \mathbb{P}^1\setminus B$ with $j=1,2$, let
$\{x_{j1},\ldots,x_{jn}\}$ be a labelling of the fiber above~$q_j$
inducing a bijection $\nu_j$. Then the bijections $\nu_1$ and
$\nu_2$ are said to be equivalent if there exists a path
$\mu\subset\mathbb{P}^1\setminus B$ running from $q_1$ to $q_2$
such that the lift of $\mu$ to $C$ with initial point $x_{1i}$ has
end point $x_{2i}$, for $i=1,\ldots,n$. The equivalence class
$[\nu]$ of an induced bijection $\nu$ is called a labelling class
for $p$.

\begin{definition}
\label{prymtriple} Consider the triple $\mathcal{T}=(A^{\oplus
m},p,[\nu])$ for $m\geq 1$, where $p:C\rightarrow \mathbb{P}^1$ is
a covering of degree $md$ and $[\nu]$ is a labelling class for
$p$. We say that $\mathcal{T}$ is a {\it Prym triple} if ${\rm
Aut}(A^{\oplus m})$ contains the monodromy group of~$p$ with
respect to $[\nu]$.
\end{definition}

Let $B=\{b_1,\ldots,b_n\}$ be a finite subset of $\mathbb{P}^1$
and take non-trivial permutations $\sigma_1,\ldots,\sigma_n\in
{\rm Aut}(A^{\oplus m})$ such that
$\sigma_n\cdot\ldots\cdot\sigma_1=(1)$ and
$G=\langle\sigma_1,\ldots,\sigma_n\rangle$ is a transitive
subgroup of $S_{md}$. Recalling the monodromy version of RET (cf.\
\cite{Mi}, p.\ 92), we may assume that $p:C\rightarrow
\mathbb{P}^1$ is a covering of degree $md$ with branch locus $B$
and labelling class $[\nu]$ such that the ramification of $p$
above $b_i$ is induced by $\sigma_i$, for $i=1,\ldots,n$. Then
$(A^{\oplus m},p,[\nu])$ is a Prym triple.\\

\noindent \textbf{A symmetric correspondence.} Assume that
$\mathcal{T}=(A^{\oplus m},p,[\nu])$ is a Prym triple with
$p:C\rightarrow \mathbb{P}^1$ and $\nu:p^{-1}(q)\rightarrow
\{1,\ldots,md\}$. We shall define a symmetric correspondence on
$C\times C$. Take the Galois closure $\pi:X\rightarrow
\mathbb{P}^1$ of $p$ and let $H$ be the Galois group of the
covering $X\rightarrow C$. Denote $G={\rm Gal}(\pi)$ and choose a
Galois labelling $\pi^{-1}(q)\leftrightarrow G$ such that the
induced labelling $p^{-1}(q)\leftrightarrow H\backslash G$ yields
$\nu(H)=1$. Define the group $\Sigma=\{\sigma_g\in S_{md}|\,g\in
G\}$, where $\sigma_g$ is the permutation sending $\nu(Hf)\mapsto
\nu(Hfg^{-1})$ and consider the canonical homomorphism
$\phi:G\rightarrow \Sigma$, $g\mapsto \sigma_g$. Since ${\rm
ker}(\phi)$ is a normal subgroup of $G$ contained in $H$, the
minimality of $\pi$ dictates that ${\rm ker}(\phi)$ is trivial,
implying that $\phi$ is an isomorphism. Noticing that, with
respect to the Galois labelling, the monodromy group of $\pi$ acts
on $G$ via multiplication on the right with the elements of $G$,
we conclude that $\Sigma$ is the monodromy group of $p$ with
respect to $[\nu]$. Hence $\phi:G\hookrightarrow {\rm
Aut}(A^{\oplus m})$ is a transitive representation. The fact that
$H$ is the stabilizer of the letter 1 with respect to $\phi$ thus
implies that $\mathcal{P}=(A^{\oplus m},r,\pi,\phi)$ for $r\in
\{r_{\!{}_+},r_{\!{}_-}\}$ represents Prym data. Therefore
$$
D_{\!{}_\mathcal{T}}:=\hat{D}_{\!{}_\mathcal{P}}+r\Delta_C
$$ is a well-defined symmetric correspondence on $C\times C$.
Recall that $k$ is the eigenvalue of the eigenvector
$(1,\ldots,1)$ of $A$. Because $\nu(Hg)=(\phi(g))^{-1}(1)$ for all
$g\in G$, equation \eqref{matrixidentity} gives the following
interpretation of $D_{\!{}_\mathcal{T}}$.

\begin{lemma}
\label{matrixcorrespondence} Let $\{x_1,\ldots,x_{md}\}$ be a
labelling in the class $[\nu]$ and denote $A^{\oplus
m}=(s_{ij})_{i,j=1}^{md}$. Then the point $(x_i,x_j)$ appears in
$D_{\!{}_\mathcal{T}}$ with multiplicity~$s_{ij}$. In particular,
$D_{\!{}_\mathcal{T}}$ is of bidegree $(k,k)$. \qed
\end{lemma}

In fact, if $S$ denotes the set of non-zero entries of $A$ and for
each $s\in S$ we define a set $\hat{G}_s=\{\hat{g}\in H\backslash
G/H|\,s_{(\phi(g))(1),1}=s\,\}$, then we find reduced divisors
$D_s=\sum_{\hat{g}\in \hat{G}_s}\hat{\Gamma}_g$ on $C\times C$
such that $D_{\!{}_\mathcal{T}}=\sum_{s\in S}s\,D_s$.

\begin{definition}
\label{prym2} Let
$\gamma_{\!{}_\mathcal{T}}=\hat{\gamma}_{\!{}_\mathcal{P}}+ r{\rm
id}_{J(C)}$. Then $P_{\!{}_\pm}(\mathcal{T})={\rm ker}
(\gamma_{\!{}_\mathcal{T}}-r_{\!{}_\mp}{\rm id}_{J(C)})_{\!{}_0}$
are the Prym varieties associated to $\mathcal{T}$.
\end{definition}

\noindent \textbf{Remark.} Given a second labelling $\nu'$ of the
fiber $p^{-1}(q)$, according to RET there exists an $f\in {\rm
Aut}(p)$ such that $\nu=\nu'\circ f$ if and only if $\nu$ and
$\nu'$ yield the same monodromy representation for $p$. If such an
$f$ exists, then $f$ induces an isomorphism of the Prym
varieties.\\

By Proposition \ref{complement}, if $m=1$, then
$P_{\!{}_+}(\mathcal{T})$ and $P_{\!{}_-}(\mathcal{T})$ are
complementary subvarieties of $J(C)$ defined by
$P_{\!{}_\pm}(\mathcal{T})={\rm
im}(\gamma_{\!{}_\mathcal{T}}-r_{\!{}_\pm}{\rm id}_{J(C)})$. If in
addition $D_{\!{}_\mathcal{T}}$ is fixed point free and
$r_{\!{}_+}=1$, then a theorem of Kanev (cf.\ \cite{Ka1}, Theorem
3.1) states that $P_{\!{}_+}(\mathcal{T})$ is a Prym-Tyurin
variety of exponent $1-r_{\!{}_-}$ for $C$.\\

We now try to compute the dimension of
$P_{\!{}_\pm}(\mathcal{T})$. Let $\eta\in \{0,1\}$ be such that
$\eta=1$ if $k\neq r_{\!{}_\pm}$ and $\eta=0$ else.

\begin{proposition}
\label{dimension} Let $\mathcal{T}=(A^{\oplus m},p,[\nu])$ be a
Prym triple with $p:C\rightarrow \mathbb{P}^1$ and assume that $A$
has diagonal $(s,\ldots,s)$. Denote $\mathcal{T}'=(A^{\oplus
m}-sI_{md},p,[\nu])$ and $\mathcal{T}_0=(J_d^{\oplus m},p,[\nu])$.
Using the notation $d_{\!{}_\pm}={\rm
dim}\,P_{\!{}_\pm}(\mathcal{T})$ and $d_0={\rm
dim}\,P_{\!{}_-}(\mathcal{T}_0)$, we have the following identity
for the dimension of $P_{\!{}_\pm}(\mathcal{T})$:
$$
\pm(r_{\!{}_+}-r_{\!{}_-})d_{\!{}_\pm}=(k-r_{\!{}_\pm})\eta d_0+
(r_{\!{}_\pm}-s)g(C)-k+s
+\frac{1}{2}(D_{\mathcal{T}'}\cdot\Delta_C)\,.
$$
\end{proposition}

\proof Denote $\gamma_{{}_0}=\gamma_{\!{}_{\mathcal{T}_0}}$ and
define an endomorphism $\varepsilon$ on $J(C)$ such that
$\varepsilon=\gamma_{\!{}_\mathcal{T}}-k\,{\rm id}_{J(C)}$ if
$\eta=1$ and $\varepsilon={\rm id}_{J(C)}$ else. Then Lemma
\ref{matrixcorrespondence} implies that
$\varepsilon(\gamma_{\!{}_\mathcal{T}}-r_{\!{}_+}{\rm
id}_{J(C)})(\gamma_{\!{}_\mathcal{T}}-r_{\!{}_-}{\rm
id}_{J(C)})=0$. Recall that $\varrho_a(\gamma_{\!{}_\mathcal{T}})$
and $\varrho_a(\gamma_{{}_0})$ are self-adjoint w.r.t.\ the form
$c_1(\Theta_C)$. Hence by direct consequence of Proposition
\ref{isogeny}, if $\eta=1$ (resp.\ $\eta=0$), then
$\varrho_a(\gamma_{\!{}_\mathcal{T}})$ has eigenvalues
$k,r_{\!{}_+},r_{\!{}_-}$ (resp.\ $r_{\!{}_+},r_{\!{}_-}$) with
respective multiplicities $d_0,d_{\!{}_-},d_{\!{}_+}$ (resp.\
$d_{\!{}_-},d_{\!{}_+}$). Since ${\rm
Tr}\bigl(\varrho_a(\gamma_{\!{}_\mathcal{T}})\bigr)={\rm
Tr}\bigl(\varrho_a(\gamma_{\!{}_{\mathcal{T}'}})\bigr)+ s\,g(C)$,
we obtain
\begin{equation}
\label{pair} \left\{ \begin{array}{l}
\,\,\,\,\,\,d_{\!{}_+}+\,\,\,\,\,d_{\!{}_-}+\,\,\,\,\eta d_0=g(C) \\
r_{\!{}_-}d_{\!{}_+}+r_{\!{}_+}d_{\!{}_-}+k\eta d_0={\rm
Tr}\bigl(\varrho_a(\gamma_{\!{}_{\mathcal{T}'}})\bigr)+s\,g(C)
\end{array}\right.
\end{equation} Let
${\rm Tr}_r(\gamma_{\!{}_{\mathcal{T}'}})$ be the rational trace
of $\gamma_{\!{}_{\mathcal{T}'}}$, i.e., ${\rm
Tr}_r(\gamma_{\!{}_{\mathcal{T}'}})$ is the trace of the extended
rational representation $(\varrho_r\otimes 1)
(\gamma_{\!{}_{\mathcal{T}'}})$ of $H^1(C,\mathbb{Z})\otimes
\mathbb{C}$. As $\varrho_r\otimes 1$ is equivalent to
$\varrho_a\oplus \overline{\varrho_a}$, it follows that ${\rm
Tr}_r(\gamma_{\!{}_{\mathcal{T}'}})=2\,{\rm
Tr}\bigl(\varrho_a(\gamma_{\!{}_{\mathcal{T}'}})\bigr)$. With
$D_{\mathcal{T}'}$ being of bidegree $(k-s,k-s)$, Proposition
11.5.2. of \cite{BL} implies
$$
{\rm
Tr}\bigl(\varrho_a(\gamma_{\!{}_{\mathcal{T}'}})\bigr)=\frac{1}{2}
{\rm Tr}_r(\gamma_{\!{}_{\mathcal{T}'}})=
k-s-\frac{1}{2}(D_{\mathcal{T}'}\cdot \Delta_C)\,.
$$ Solving \eqref{pair} for $d_{\!{}_\pm}$ we obtain the desired
result.\qed\\

Proposition \ref{complement} implies that for $m=1$ we have
$d_0=0$. To compute ${\rm dim}\,P_{\!{}_\pm}(\mathcal{T})$ we need
to determine the intersection number
$(D_{\mathcal{T}'}\cdot\Delta_C)$. We shall do this for a Prym
triple $\mathcal{T}=(A^{\oplus m},p,[\nu])$, where $A^{\oplus
m}=(s_{i,j})_{i,j=1}^{md}$ has zero diagonal. Denoting the set of
nonzero entries of $A$ by $S$, we recall that
$D_{\!{}_\mathcal{T}}=\sum_{s\in S}s\,D_s$ with $D_s$ reduced.
Hence it suffices to determine the local intersection numbers
$(D_s\cdot\Delta_C)_{(x,x)}$ at $(x,x)$ for a ramification point
$x\in C$ of $p:C\rightarrow \mathbb{P}^1$. Let $b\in \mathbb{P}^1$
be the corresponding branch point and assume that the local
monodromy at $b$ is given by $\sigma_b\in S_{md}$. Further, let
$\tau\in S_{md}$ be the cycle factor of $\sigma_b$ which describes
the ramification at $x$ and assume that it is of order $l$. For
each $s\in S$ we define a set $T_{\tau,s}$ of elements $t\in
\{1,\ldots,l-1\}$ for which there exists a $j\in \{1,\ldots,md\}$
such that $s_{j,\tau^t(j)}=s$. Then:

\begin{lemma}
\label{intersectionnumber} For each $s\in S$ we have
$(D_s\cdot\Delta_C)_{(x,x)}=\#(T_{\tau,s})$.
\end{lemma}

\proof After a suitable choice of coordinates on a small open
neighborhood of $x$ the covering $p$ is given by $z\mapsto z^l$.
Then near the point $(x,x)$ the reduced divisor $D_s$ can be
described as the union of graphs of the multiplications $z\mapsto
\zeta_l^t z$ (for $t\in T_{\tau,s}$) with $\zeta_l={\rm
exp}(\frac{2\pi}{l}\sqrt{-1})$. Obviously these graphs intersect
$\Delta_C$ transversally in $(x,x)$, thus implying
$(D_s\cdot\Delta_C)_{(x,x)}=\#(T_{\tau,s})$.\qed\\

Hence, given the branch locus $B$ of $p$ and, for each $b\in B$,
the set $R_b$ of cycle factors in the cycle decomposition of
$\sigma_b$, we can calculate the intersection number as a sum
$$
(D_{\!{}_\mathcal{T}}\cdot \Delta_C)=\sum_{b\in B}\sum_{\tau\in
R_b}\sum_{s\in S}s\,\#(T_{\tau,s})\,.
$$ For an application of the lemma we refer to Example
\ref{twisted}.

\end{section}

\begin{section}{Examples}
\label{examples}

Our first example has been covered by Kanev in \cite{Ka2}. We will
treat it by a different method.

\begin{example}
\label{schlaefli} \textbf{(Schl\"{a}fli graph)} Let $\mathcal{L}$
be the intersection graph of the 27 lines on a non-singular cubic
surface in $\mathbb{P}^3$. In the notation of Example
\ref{schlaefligraph} we take $\tau_i$ $(i=1,\ldots,5)$ (resp.\
$\tau_6$) to be the transformation that interchanges the rows of
the double-six $M_{i,i+1}$ (resp.\ $M_{1,2,3}$). Denoting the
adjacency matrix of $\mathcal{L}$ by $A$, we recall that ${\rm
Aut}(A)=\langle\tau_1,\ldots,\tau_6\rangle$ is a transitive
subgroup of $S_{27}$. Note moreover that each $\tau_j$ is a
reflection with exactly 15 fixed points.

Let $n\geq 7$ be an integer and choose a subset
$B=\{b_1,\ldots,b_{2n}\}$ of $\mathbb{P}^1$. Then we know that
there exists a Prym triple $\mathcal{T}=(A,p,[\nu])$ for a
covering $p:C\rightarrow \mathbb{P}^1$ with branch locus $B$ and
monodromy group ${\rm Aut}(A)$ such that its ramification over
$b_i$ is induced by a $\tau_j$. According to Hurwitz' formula the
curve $C$ is of genus $6n-26$. Since no vertex $v$ of
$\mathcal{L}$ is adjacent to $\tau_j(v)$, it follows that
$D_{\!{}_\mathcal{T}}$ is fixed point free. As $A$ has eigenvalues
$k=10$, $r_{\!{}_+}=1$ and $r_{\!{}_-}=-5$, Proposition
\ref{dimension} implies that $P_{\!{}_+}(\mathcal{T})$ is an
$(n-6)$-dimensional Prym-Tyurin variety of exponent 6 for the
curve $C$. With regard to moduli, note that for $n=12$ we have
$g(C)=46$, ${\rm dim}\,P_{\!{}_+}(\mathcal{T})=6$ and $\#(B)={\rm
dim}\,\mathcal{A}_6+{\rm dim\,Aut}(\mathbb{P}^1)=24$, where
$\mathcal{A}_6$ is the moduli space of 6-dimensional principally
polarized abelian varieties.
\end{example}

\begin{example}
\label{lattice} \textbf{(Lattice graphs)} For $n\geq 3$ let $A$ be
the adjacency matrix of the lattice graph $L_2(n)$ with vertex set
$\{1,\ldots,n\}^2$. We define a group
$G=\langle\varphi_0,\varphi_1,\varphi_2,\varphi_3\rangle$
generated by reflections $\varphi_h:=(\tau_h,\tau_h^{-1})\circ t$
in ${\rm Aut}(A)$, where $t$ acts on $\{1,\ldots,n\}^2$ by
exchange of coordinates and $\tau_0,\tau_1,\tau_2,\tau_3\in S_n$
are given by $\tau_0=(1)$, $\tau_1=(1\hskip0.25cm n)$,
$\tau_2=(2\hskip0.25cm n)$ and $\tau_3=(1\hskip0.25cm
2\hskip0.25cm \cdots\hskip0.25cm n)$. Then $G$ is a transitive
subgroup of ${\rm Aut}(A)$; indeed, identifying $\{1,\ldots,n\}$
and $\mathbb{Z}/n\mathbb{Z}$, we have
\begin{itemize}
\item[a)] $(\varphi_1\circ \varphi_3)^m(1,1)=(1,h+1)$ for
$m=1,\ldots,n-2$; \item[b)] $(\varphi_3\circ
\varphi_0)(i,j)=(i+1,j-1)$ for $i,j=1,\ldots,n$ with $j\neq n-i$;
\item[c)] $\bigl((\varphi_3\circ
\varphi_0)^{m-1}\circ\varphi_2\bigr)(2,1)=(m,n-m+1)$ for
$m=1,\ldots,n$.
\end{itemize} Given an integer $l\geq 0$, choose a subset
$B=\{\{b_1,\ldots,b_{2l+8}\}$ of $\mathbb{P}^1$. We may assume
that $\mathcal{T}=(A,p,[\nu])$ is a Prym triple for a covering
$p:C\rightarrow \mathbb{P}^1$ with branch locus $B$ and monodromy
group $G$ such that its ramification over $b_i$ is induced by a
$\varphi_h$. Then $C$ is of genus $(n-1)^2+\frac{1}{2}ln(n-1)$. As
$D_{\!{}_\mathcal{T}}$ is fixed point free and $A$ has eigenvalues
$k=2(n-1)$, $r_{\!{}_+}=n-2$ and $r_{\!{}_-}=-2$, it follows that
$$
{\rm dim}\,P_{\!{}_+}(\mathcal{T})=
(n-1)(n-3)+\frac{1}{2}l(n-1)(n-2)\,.
$$ Hence, for $n=3$ and $l\geq 1$ we obtain a finite number of finite
dimensional families of $l$-dimensional Prym-Tyurin varieties of
exponent 3 for curves of genus $3l+4$. In anticipation of section
\ref{exponent} we shall say that $\mathcal{T}$ is of {\it type}
$l$ whenever $n=3$ and $l\geq 1$.
\end{example}

\begin{example}
\label{twisted} For an example involving symmetric correspondences
with fixed points, let $n\geq 3$ and assume that $A$ is the
adjacency matrix of the graph $\overline{L_2(n)}$ (the complement
of $L_2(n)$) with vertex set $\{1,\ldots,n\}^2$. Assume that $t\in
{\rm Aut}(A)$ acts on $\{1,\ldots,n\}^2$ by exchange of
coordinates and for $h=1,\ldots,n-1$ define the reflection
$\sigma_h:=\bigl((1\hskip0.25cm h+1),(1)\bigr)$ in $S_n\times
S_n$. We observe that ${\rm Aut}(A)$ is generated by the elements
$t$ and $\sigma_1,\ldots,\sigma_{n-1}$. Clearly, no vertex
$(i,j)\in \{1,\ldots,n\}^2$ is adjacent to $\sigma_h(i,j)$.
Further, $(i,j)$ is adjacent to $t(i,j)$ if and only if $i\neq j$.

Given nonnegative integers $l_1,l_2$, we choose a subset
$B=B_1\sqcup B_2$ of $\mathbb{P}^1$ with
$B_1=\{b_{1,1},\ldots,b_{1,2(l_1+1)}\}$ and
$B_2=\{b_{2,1},\ldots,b_{2,2(l_2+n-1)}\}$. Let
$\mathcal{T}=(A,p,[\nu])$ be a Prym triple for a covering
$p:C\rightarrow \mathbb{P}^1$ with branch locus $B$ and monodromy
group ${\rm Aut}(A)$ such that its ramification over $b_{1,i}$
(resp.\ $b_{2,j}$) is induced by $t$ (resp.\ some $\sigma_h$).
Then the curve $C$ is of genus
$\frac{1}{2}(n-1)(n-2)+\frac{1}{2}l_1n(n-1)+l_2n$ and Lemma
\ref{intersectionnumber} implies that $(D_{\!{}_\mathcal{T}}\cdot
\Delta_C)=(l_1+1)(n-1)n$. It follows that
$P_{\!{}_+}(\mathcal{T})$ is of dimension $l_1(n-1)+l_2$.

In view of moduli, note that for $l_1=0$ and $n\geq
\frac{1}{4}(l_2^2-3l_2+6)$ we have ${\rm
dim}\,P_{\!{}_+}(\mathcal{T})=l_2$ and $\# B\geq {\rm
dim}\,\mathcal{A}_{l_2}+{\rm dim\,Aut}(\mathbb{P}^1)$. In
particular, if $l_1=0$ and $n=l_2=6$, then $g(C)=46$. Moreover,
since $S_2\rtimes (S_n\times S_n)$ has no subgroup of index $n$,
Galois theory implies that no factorization $p:C\xrightarrow{n:1}
C'\xrightarrow{n:1} \mathbb{P}^1$ exists.
\end{example}

\begin{example}
\label{latinsquare} \textbf{(Latin square graphs)} Given an
integer $n\geq 3$, we assume that $A$ is the adjacency matrix of
the Latin square graph $L_3(n)$. We recall from Example
\ref{latinsquaregraphs} that $(\mathbb{Z}/n\mathbb{Z})^2$ induces
a transitive subgroup of ${\rm Aut}(A)$ via translation; as such
it coincides with $\langle (1,1),(1,2)\rangle$. Viewed as
permutations of the vertex set $(\mathbb{Z}/n\mathbb{Z})^2$, the
translations $(1,1)$ and $(1,2)$ split into $n$ mutually disjoint
$n$-cycles. For $n\geq 4$ the vertices $(i+1,j+1)$ and $(i+1,j+2)$
of $L_3(n)$ are non-adjacent to $(i,j)$.

Now assume that $n\geq 4$. We choose an integer $l\geq 2$ and a
subset $B=\{b_1,\ldots,b_{ln}\}$ of $\mathbb{P}^1$. Then there
exists a Prym triple $\mathcal{T}=(A,p,[\nu])$ for a covering
$p:C\rightarrow \mathbb{P}^1$ with branch locus $B$ and monodromy
group $(\mathbb{Z}/n\mathbb{Z})^2$ such that its ramification over
$b_i$ is induced by $(1,1)$ or $(1,2)$. We find that $C$ is of
genus $1-n^2+\frac{1}{2}ln^2(n-1)$. Moreover, since ${\rm
deg}(p)=\# (\mathbb{Z}/n\mathbb{Z})^2$, it is immediately seen
that $p$ is a Galois covering. Using the fact that
$D_{\!{}_\mathcal{T}}$ is fixed point free and $A$ has eigenvalues
$k=3(n-1)$, $r_{\!{}_+}=n-3$ and $r_{\!{}_-}=-3$, we compute
$$
{\rm
dim}\,P_{\!{}_+}(\mathcal{T})=-(n-1)(n-2)+\frac{1}{2}ln(n-1)(n-3)\,.
$$ Hence, for $n=4$ we get finitely many finite dimensional
families of $6(l-1)$-dimensional Prym-Tyurin varieties of exponent
4 for curves of genus $24l-15$.
\end{example}

\end{section}

\begin{section}{A splitting}
\label{splitting}

We show that for certain Prym triples $\mathcal{T}=(A^{\oplus
m},p,[\nu])$ the covering $p:C\rightarrow \mathbb{P}^1$ splits
into a covering $f:C\rightarrow C'$ of degree $d$ and a covering
$h:C'\rightarrow \mathbb{P}^1$ of degree $m$ such that $f$ depends
essentially on $D_{\!{}_\mathcal{T}}$. Recall that $k$ is the
eigenvalue of the eigenvector $(1,\ldots,1)$ of $A$.

\begin{theorem}
\label{factorization} Assume that $A\in \{0,1\}^{d\times d}$ is a
Prym matrix with zero diagonal and eigenvalue $k$ of multiplicity
1. Given $m\geq 2$, let $\mathcal{T}=(A^{\oplus m},p,[\nu])$ be a
Prym triple for a covering $p:C\rightarrow \mathbb{P}^1$ with
branch locus $B$. Denote $C_0:=p^{-1}(\mathbb{P}^1\setminus B)$.
Then there exists a unique splitting
$$
p:C\xrightarrow[f]{d:1} C'\xrightarrow[h]{m:1} \mathbb{P}^1
$$ such that, for any $(x,x')\in (C_0\times C_0)\setminus\Delta_{C_0}$,
the points $x,x'$ are in the same fiber of $f$ if and only if
there is a finite sequence of points $x=x_0,\ldots,x_l=x'$ on
$C_0$ with $(x_j,x_{j+1})\in D_{\!{}_\mathcal{T}}$ for all
$j=0,\ldots,l-1$.
\end{theorem}

\proof Fix a point $q_0\in \mathbb{P}^1\setminus B$ and assume
that $\nu$ is a labelling of the fiber $p^{-1}(q_0)$. We denote
$S=\{1,\ldots,m\}$, $T=\{1,\ldots,d\}$ and identify $S\times T$
with $\{1,\ldots,md\}$ via the bijection $(s,t)\leftrightarrow
(s-1)m+t$. Then $\nu$ turns into a bijection
$(\nu_1,\nu_2):p^{-1}(q_0)\rightarrow S\times T$. Let $\Sigma$ be
the monodromy group of $p$ with respect to $(\nu_1,\nu_2)$ and
split its elements accordingly into $\sigma=(\sigma_1,\sigma_2)$.
Denoting $A:=(a_{i,j})_{i,j=1}^d$, we may view $A^{\oplus m}$ as
the matrix of entries $c_{u,u'}$, where $u=(s,t)$ and $u'=(s',t')$
run through the set $S\times T$, such that $c_{u,u'}=a_{t,t'}$ if
$s=s'$ and $c_{u,u'}=0$ else. According to Proposition \ref{srg}
the matrix $A$ is the adjacency matrix of a connected strongly
regular graph $\mathcal{G}$ on $d$ vertices. Thus, for $u=(s,t)$
and $u'=(s',t')$ there exists a finite sequence
$u=u_0,\ldots,u_l=u'$ in $S\times T$ such that $c_{u_j,u_{j+1}}=1$
for all $j=0,\ldots,l-1$ if and only if $s=s'$. Hence
$\sigma_1(\,\cdot\,,t)=\sigma_1(\,\cdot\,,t')$ for all $\sigma\in
\Sigma$, i.e., there is a unique $\tau_\sigma\in S_m$ such that
$\sigma_1(\,\cdot\,,t)=\tau_\sigma(\,\cdot\,)$ for
all $t\in T$.\\

Let $\pi:X\rightarrow \mathbb{P}^1$ be the Galois closure of $p$
and denote $G={\rm Gal}(\pi)$. As we have seen in section
\ref{nongalois}, there exists an isomorphism $\phi:G\rightarrow
\Sigma$ such that the Galois group $H$ of $X\rightarrow C$ is the
stabilizer of $(1,1)\in S\times T$ w.r.t.\ $\phi$ and any Galois
labelling of a fiber of $\pi$ induces a labelling in the class
$[\nu]$ via the identification $Hg\leftrightarrow g^{-1}(1,1)$. We
let $H'$ be the stabilizer of $1\in S$ with respect to
$\psi\circ\phi$, where $\psi:\Sigma\rightarrow S_m$ is the
transitive representation induced by $\sigma\mapsto \tau_\sigma$.
Write $C'=X/H'$; since $H\subset H'$ (resp.\ $H'\subset G$) is a
subgroup of index $d$ (resp.\ $m$), there are canonical coverings
$f:C\rightarrow C'$ of degree $d$ and $h:C'\rightarrow
\mathbb{P}^1$ of degree $m$ such that $p=h\circ f$. Take a point
$q\in \mathbb{P}^1\setminus B$ and a Galois labelling
$\pi^{-1}(q)\leftrightarrow G$. For any element $g\in G$, if
$(s,t)=g^{-1}(1,1)$, then $H'g\leftrightarrow g^{-1}(1)=s$, i.e.,
on the fiber $p^{-1}(q)$ the covering $f$ is given by
$(s,t)\mapsto s$. With reference to Lemma
\ref{matrixcorrespondence} we conclude that $f$ has the desired
properties. Using the monodromy of $p$, the reader will easily
check that the splitting is unique.\qed\\

With $\mathcal{T}$, $f$ and $h$ as above, we say that the pair of
coverings $(f,h)$ represents the {\it canonical splitting} for
$\mathcal{T}$.

\begin{corollary}
\label{usualprym1} For integers $d,m\geq 2$, assume that
$\mathcal{T}=\bigl((J_d-I_d)^{\oplus m},p,[\nu]\bigr)$ is a Prym
triple associated to a covering $p:C\rightarrow \mathbb{P}^1$ and
let $(f,h)$ be its canonical splitting. Then
$P_{\!{}_+}(\mathcal{T})$ is the usual Prym variety associated to
the covering $f$, i.e., $P_{\!{}_+}(\mathcal{T})$ and ${\rm
im}f^*$ are complementary subvarieties of $J(C)$.
\end{corollary}

\proof According to Theorem \ref{factorization} we have
$D_{\!{}_\mathcal{T}}(x)=-x+f^*f(x)$ for all $x\in C$ in an
unramified fiber of $p$. Hence $\gamma_{\!{}_\mathcal{T}}+{\rm
id}_{J(C)}=f^*N_f$ and thus $P_{\!{}_-}(\mathcal{T})={\rm
im}\,f^*$. As $(J_d-I_d)^{\oplus m}$ is a Prym matrix, Proposition
\ref{complement} implies that $P_{\!{}_+}(\mathcal{T})$ and ${\rm
im}f^*$ are complementary in $J(C)$.\qed\\

Corollary \ref{usualprym1} has a natural converse. Before
addressing this, we recall that a smooth projective curve of genus
$g$ is $m$-gonal for all $m\geq \lceil \frac{g}{2} \rceil +1$
(cf.\ \cite{ACGH}, Existence Theorem, p.\ 206).

\begin{corollary}
\label{usualprym2} Assume that $f:C\rightarrow C'$ is a covering
of degree $d\geq 2$ of a curve $C'$ of genus $g\geq 1$ and let
$h:C'\rightarrow \mathbb{P}^1$ be a covering of degree $m\geq
\lceil \frac{g}{2} \rceil +1$. Then there exists a labelling class
$[\nu]$ for the covering $h\circ f$ such that
$\mathcal{T}=\bigl((J_d-I_d)^{\oplus m},h\circ f,[\nu]\bigr)$ is a
Prym triple and $P_{\!{}_+}(\mathcal{T})$ is the usual Prym
variety associated to $f$.
\end{corollary}

\proof Take a point $q\in \mathbb{P}^1$ outside the branch locus
of $h\circ f$. Then we can define a bijection
$\nu=(\nu_1,\nu_2):(h\circ f)^{-1}(q)\rightarrow
\{1,\ldots,m\}\times \{1,\ldots,d\}$ such that
$\nu_1(x)=\nu_1(x')\Longleftrightarrow f(x)=f(x')$, for all
$x,x'\in (h\circ f)^{-1}(q)$. It is immediately seen that
$\bigl((J_d-I_d)^{\oplus m},h\circ f,[\nu]\bigr)$ represents a
Prym triple with canonical splitting $(f,h)$. Now apply Corollary
\ref{usualprym1}.\qed \\

Given integers $d,m\geq 2$, assume that
$\mathcal{T}=\bigl((J_d-I_d)^{\oplus m},p,[\nu]\bigr)$ is a Prym
triple. We shall call $\mathcal{T}$ {\it simple} if its canonical
splitting $(f,h)$ is simple, i.e., if $f$ and $h$ are simply
branched coverings such that no ramified fiber of $h$ contains a
branch point of $f$ and no unramified fiber of $h$ contains more
than one branch point of $f$. It should be noted that simplicity
can also be described in terms of the monodromy of $p$ alone,
without reference to $f$ and~$h$.\\

To conclude this section, we use simple Prym triples to
characterize (at least up to isogeny) abelian varieties
corresponding to the general popints of $\mathcal{A}_4$ and
$\mathcal{A}_5$.

\begin{lemma}
\label{A4A5} (1) The general 4-dimensional principally polarized
abelian variety is isogenous to a Prym variety
$P_{\!{}_+}(\mathcal{T})$ for a simple Prym triple
$\mathcal{T}=\bigl((J_2-I_2)^{\oplus 3},p,[\nu]\bigr)$ such that
the covering $p$ has exactly 4 simple and 10 double branch points.

\noindent (2) The general 5-dimensional principally polarized
abelian variety is of the form $P_{\!{}_+}(\mathcal{T})$, where
$\mathcal{T}=\bigl((J_2-I_2)^{\oplus 4},p,[\nu]\bigr)$ is a simple
Prym triple such that the covering $p$ has exactly 18 double
branch points.
\end{lemma}

\proof For integers $g\geq 1$ and $n\geq 0$, let
$\mathcal{R}(g,n)$ be the moduli space of equivalence classes of
double coverings $f:C\rightarrow C'$ with $C'$ of genus $g$ and
$f$ branched at $2n$ distinct points of $C'$. We shall need the
following fact: Let $m$ be an integer. If $m\geq \lceil
\frac{g}{2} \rceil +1$, then for a double covering $f:C\rightarrow
C'$ corresponding to a general point of $\mathcal{R}(g,n)$ there
exists an $m$-fold covering $h:C'\rightarrow \mathbb{P}^1$ such
that the covering pair $(f,h)$ is simple. The proof is left to the
reader. As in \cite{BCV}, p.\ 122, we let
$p_{(g,n)}:\mathcal{R}(g,n)\rightarrow
\mathcal{A}_{g+n-1}(\delta)$ be the usual Prym morphism, where
$\mathcal{A}_{g+n-1}(\delta)$ is the moduli space of abelian
$g$-folds with polarization type $\delta$. According to
\cite{BCV}, Theorem 2.2, the morphism
$p_{(3,2)}:\mathcal{R}(3,2)\rightarrow \mathcal{A}_4(1,2,2,2)$ is
dominant. Moreover, for the general double covering
$f:C\rightarrow C'$ with 4 branch points and $g(C')=3$ there
exists a 3-fold covering $h:C'\rightarrow \mathbb{P}^1$ such that
the pair $(f,h)$ is simple and the covering $h\circ f$ has exactly
4 simple and 10 double branch points. Together with Corollary
\ref{usualprym2} this shows (1). To prove (2) we recall that
$p_{(6,0)}:\mathcal{R}(6,0)\rightarrow \mathcal{A}_5$ is dominant
(cf.\ \cite{Mu}). Hence it suffices to note that for the general
\'{e}tale double covering $f:C\rightarrow C'$ with $g(C')=6$ there
exists a 4-fold covering $h:C'\rightarrow \mathbb{P}^1$ branched
at 18 points such that the pair $(f,h)$ is simple. \qed

\end{section}

\begin{section}{Prym-Tyurin varieties of arbitrary exponent $\geq 3$}
\label{exponent}

We show how the graph $\overline{L_2(n)}\in {\rm
SRG}\bigl(n^2,(n-1)^2,(n-2)^2,(n-1)(n-2)\bigr)$ for $n\geq 3$ can
be employed to construct families of Prym-Tyurin varieties of
exponent $n$. These varieties turn out to be the product of the
Jacobians of two $n$-gonal curves.

\begin{example}
\label{prymtyurin} Given an integer $n\geq 3$, we shall try to
construct Prym-Tyurin varieties of exponent $n$. Assume that $A$
is the adjacency matrix of the graph $\overline{L_2(n)}$ with
vertex set $\{1,\ldots,n\}^2$. Recall that $S_n\times S_n$ is a
transitive subgroup of ${\rm Aut}(A)$. For $i=1,\ldots,n-1$ we
define reflections $\sigma_{1,i}:=\bigl((1\hskip0.25cm
i+1),(1)\bigr)$ and $\sigma_{2,i}:=\bigl((1),(1\hskip0.25cm
i+1)\bigr)$ in $S_n\times S_n$. Note that $S_n\times S_n$ is
freely generated by these reflections.

Given nonnegative integers $l_1,l_2$ such that $l_1+l_2\geq 1$,
let $B=B_1\sqcup B_2$ be a finite subset of $\mathbb{P}^1$ with
$B_m=\{b_{m,1},\ldots,b_{m,2(l_m+n-1)}\}$. Assume that
$\mathcal{T}=(A,p,[\nu])$ is a Prym triple for a covering
$p:C\rightarrow \mathbb{P}^1$ with branch locus $B$ and monodromy
group $S_n\times S_n$ such that its ramification over $b_{m,i}$ is
induced by a $\sigma_{m,h}$; in this situation we say that
$\mathcal{T}$ is of {\it type} $(l_1,l_2)$. Since no vertex
$(i,j)$ of $\overline{L_2(n)}$ is adjacent to $\sigma_{m,h}(i,j)$,
the correspondence $D_{\!{}_\mathcal{T}}$ is fixed point free.
Moreover, it is easily seen that all $\sigma_{m,h}$ decompose into
$n$ mutually disjoint transpositions on the set
$\{1,\ldots,n\}^2$. Hence, as $A$ has eigenvalues $k=(n-1)^2$,
$r_{\!{}_+}=1$ and $r_{\!{}_-}=-n+1$, it follows that
$P_{\!{}_+}(\mathcal{T})$ is an $(l_1+l_2)$-dimensional
Prym-Tyurin variety of exponent $n$ for the curve $C$ of genus
$(n-1)^2+(l_1+l_2)n$.
\end{example}

Recall from Example \ref{lattice} that a Prym triple $\mathcal{T}$
of type $l$ yields an $l$-dimensional Prym-Tyurin variety
$P_{\!{}_+}(\mathcal{T})$ of exponent 3. We will show that for
$n=3$ the Prym-Tyurin varieties of the preceding example are the
same as those of Example \ref{lattice}. More precisely, let $A$ be
the adjacency matrix of the lattice graph $L_2(3)$ and assume that
$[\nu]$ is a labelling class for a covering $p:C\rightarrow
\mathbb{P}^1$ of degree 9. Given the isomorphism of graphs
$\xi:L_2(3) \overset \sim \rightarrow \overline{L_2(3)}$ induced
by the matrix $\bigl(
\begin{smallmatrix} 1 & -1 \\ 1 & 1 \end{smallmatrix} \bigr)$
as in Example \ref{latticegraphs}, we have:

\begin{lemma}
\label{types} Let $l\geq 1$ be an integer. Then
$\mathcal{T}=(A,p,[\nu])$ is a Prym triple of type $l$ if and only
if there exist integers $l_1,l_2\geq 0$ such that $l_1+l_2=l$ and
$\mathcal{T}'=(J_9-I_9-A,p,[\xi^{-1}\circ\nu])$ is a Prym triple
of type $(l_1,l_2)$. In particular, if $\mathcal{T}$ is of type
$l$, then $P_{\!{}_\pm}(\mathcal{T})=P_{\!{}_\pm}(\mathcal{T}')$.
\end{lemma}

\proof It suffices to note that the reflections
$\varphi_0,\varphi_1,\varphi_2,\varphi_3\in {\rm
Aut}\bigl(L_2(3)\bigr)$ of Example \ref{lattice} satisfy the
identities $\xi^{-1}\sigma_{1,1}\xi=\varphi_0$,
$\xi^{-1}\sigma_{1,2}\xi=\varphi_0\varphi_3\varphi_0^{-1}$,
$\xi^{-1}\sigma_{2,1}\xi=\varphi_1\varphi_2\varphi_1^{-1}$ and
$\xi^{-1}\sigma_{2,2}\xi=\varphi_2$.\qed\\

For Prym-Tyurin varieties $P_{\!{}_+}(\mathcal{T})$ with
$\mathcal{T}$ of type $(l_1,l_2)$ we have the following
characterization.

\begin{theorem}
\label{characterization} Assume that $A$ is the adjacency matrix
of $\overline{L_2(n)}$ with $n\geq 3$. For nonnegative integers
$l_1,l_2$ such that $l_1+l_2\geq 1$, let $\mathcal{T}$ be a Prym
triple of type $(l_1,l_2)$ associated to $A$ and a covering
$C\rightarrow \mathbb{P}^1$. Then there exist $n$-gonal curves
$C_1$ of genus $l_1$ and $C_2$ of genus $l_2$ such that
$C=C_1\times_{{}_{\mathbb{P}^1}} C_2$ and
$P_{\!{}_+}(\mathcal{T})\simeq J(C_1)\times J(C_2)$.
\end{theorem}

\proof Let $p:C\rightarrow \mathbb{P}^1$ and $[\nu]$ be the
covering and labelling class such that $\mathcal{T}=(A,p,[\nu])$.
Then $S_n\times S_n$ is the monodromy group of $p$. We take the
inclusion $\iota:S_n\times S_n\hookrightarrow {\rm Perm}(N\times
N)$ and write $N=\{1,\ldots,n\}$. Let $\pi:X\rightarrow
\mathbb{P}^1$ be the Galois closure of $p$ and denote $G={\rm
Gal}(\pi)$. Then there exists an isomorphism $\phi:G\rightarrow
\Sigma$ such that the Galois group $H$ of $X\rightarrow C$ is the
stabilizer of $(1,1)\in N\times N$ w.r.t.\ $\iota\circ\phi$ and
any Galois labelling of a fiber of $\pi$ induces a labelling in
the class $[\nu]$ via the identification $Hg\leftrightarrow
g^{-1}(1,1)$. Take the projection mappings ${\rm pr}_1,{\rm
pr}_2:S_n\times S_n\rightarrow S_n$ onto the first and second
factor and let $H_1$ (resp.\ $H_2$) be the stabilizer of the
letter $1\in N$ w.r.t.\ $\phi_1:={\rm pr}_1\circ\phi$ (resp.\
$\phi_2:={\rm pr}_2\circ\phi$). Observing that $H=H_1\cap H_2$, we
take the quotient curves $C_m=X/H_m$ for $m=1,2$ and let
$f_m:C\rightarrow C_m$ and $h_m:C_m\rightarrow \mathbb{P}^1$ be
the canonical coverings. The transitivity of $\phi_m$ implies that
$f_m$ and $h_m$ are of degree $n$. In addition to $H=H_1\cap H_2$
we have $G=\langle H_1,H_2\rangle$. Using elementary Galois theory
we thus find
$$
\mathbb{C}(C)=\mathbb{C}(C_1)\otimes_{\mathbb{C}(\mathbb{P}^1)}
\mathbb{C}(C_2)\,.
$$ Hence $C$ is the fiber product of the $n$-gonal curves $C_1$ and
$C_2$ with projection morphisms $f_1$ and $f_2$. Let $B$ be the
branch locus of the covering $p$. Then $h_m$ for $m=1,2$ is a
simple covering with branch locus $B_m$, where $B_m$ is the set of
points $b\in B$ such that, in the notation of Example
\ref{prymtyurin}, the local monodromy of $p$ is given by a
permutation $\sigma_{m,i}$. As $\#(B_m)=2(l_m+n-1)$, we find
$g(C_m)=l_m$.\\

It remains to show that $P_{\!{}_+}(\mathcal{T})\simeq
J(C_1)\times J(C_2)$. Choose a point $q\in \mathbb{P}^1\setminus
B$ and a Galois labelling $\pi^{-1}(q)\leftrightarrow G$. Then
take a labelling $\{y_{1,1},\ldots,y_{n,n}\}$ for $p^{-1}(q)$ and
a labelling $\{z_{m,1},\ldots,z_{m,n}\}$ for $h_m^{-1}(q)$,
$m=1,2$ such that $y_{g^{-1}(1,1)}\leftrightarrow Hg$ and
$z_{m,(\phi_m(g))^{-1}(1)}\leftrightarrow H_mg$, for all $g\in G$.
We observe that $f_1^{-1}(z_{1,s})=\{\,y_{s,j}\,|\,j\in N\}$ and
$f_2^{-1}(z_{2,t})=\{\,y_{i,t}\,|\,i\in N\}$, for all $s,t\in N$.
According to Lemma \ref{matrixcorrespondence} we have
$D_{\!{}_\mathcal{T}}(y_{s,t})=\sum_{i\neq s,j\neq t}y_{i,j}$ and
therefore
$$
f_1^*z_{1,s}+f_2^*z_{2,t}=\sum_{j\in N}y_{s,j}+\sum_{i\in
N}y_{i,t} =p^*p(y_{s,t})+y_{s,t}-D_{\!{}_\mathcal{T}}(y_{s,t})\,.
$$ Hence, for $y,y'\in p^{-1}(\mathbb{P}^1\setminus B)$ and
$z_m=f_m(y)$, $z'_m=f_m(y')$ with $m=1,2$ we obtain, using divisor
class notation,
$$
f_1^*[z_1-z'_1]+ f_2^*[z_2-z'_2]=-(\gamma_{\!{}_\mathcal{T}}-{\rm
id}_{J(C)})([y-y'])\,.
$$ Consequently, defining
$\varphi=f_1^*\psi_1+f_2^*\psi_2:J(C_1)\times J(C_2)\rightarrow
J(C)$, where $\psi_m:J(C_1)\times J(C_2)\rightarrow J(C_m)$ is the
projection on the $m$-th factor, we get
$P_{\!{}_+}(\mathcal{T})={\rm im}(\gamma_{\!{}_\mathcal{T}}-{\rm
id}_{J(C)})\subset {\rm im}\,\varphi$. Because ${\rm
dim}\,P_{\!{}_+}(\mathcal{T})=l_1+l_2={\rm dim}\,J(C_1)+{\rm
dim}\,J(C_2)$, it thus follows that $\varphi:J(C_1)\times
J(C_2)\rightarrow P_{\!{}_+}(\mathcal{T})$ is an isogeny. As
$P_{\!{}_+}(\mathcal{T})$ is a Prym-Tyurin variety of exponent $n$
for $C$, the restriction of the canonical polarization $\Theta_C$
to $P_{\!{}_+}(\mathcal{T})$ is of type $(n,\ldots,n)$. Lemma
12.3.1 of \cite{BL} implies that the polarization
$\varphi^*\Theta_C$ of $J(C_1)\times J(C_2)$ is of type
$(n,\ldots,n)$, as well. Hence $\varphi:J(C_1)\times
J(C_2)\rightarrow P_{\!{}_+}(\mathcal{T})$ is an isogeny of degree
1, i.e., an isomorphism.\qed\\

\noindent \textbf{Remark.} In spite of the similarities between
the Examples \ref{latinsquare} and \ref{prymtyurin}, the preceding
theorem does not fully extend to Prym triples $\mathcal{T}$ such
as in Example \ref{latinsquare}. In fact, defining $n$-gonal
curves $C_1$ and $C_2$ analogously to those in the proof, we get
$C=C_1\times_{{}_{\mathbb{P}^1}} C_2$. A simple computation shows,
however, that the dimensions of $P_{\!{}_+}(\mathcal{T})$ and
$J(C_1)\times J(C_2)$ do not match.\\

\noindent \textbf{A different construction.} In \cite{LRR}, Lange,
Recillas and Rochas define non-trivial families of Prym-Tyurin
varieties of exponent 3. Here is a recap of their construction:
Given a hyperelliptic curve $X$ of genus $g\geq 3$ and an
\'{e}tale covering $f:\tilde{X}\rightarrow X$ of degree $3$, let
$h:X\rightarrow \mathbb{P}^1$ be the covering given by the $g_2^1$
and define the curve $C=(f^{(2)})^{-1}(g_2^1)$, where
$f^{(2)}:\tilde{X}^{(2)}\rightarrow X^{(2)}$ is the second
symmetric product of $f$. Assume for the moment that $C$ is smooth
and irreducible. Denote $\tilde{C}=\mu^{-1}(C)$, where
$\mu:\tilde{X}^2\rightarrow \tilde{X}^{(2)}$ is the canonical
$2:1$ map and let $\pi:\tilde{C}\rightarrow \tilde{X}$ be the
projection on the first factor, where $\tilde{C}$ is considered as
a curve in $\tilde{X}^{2}$. Now define the covering
$p:C\xrightarrow{9:1} \mathbb{P}^1$ induced by $h\circ f\circ
\pi:\tilde{C}\rightarrow \mathbb{P}^1$ and let $\iota:X\rightarrow
X$ be the hyperelliptic involution. Then $p$, $h\circ f$ and $h$
have the same branch locus $B$, which may be assumed to be of
cardinality $2l+8$ for $l\geq 0$. To obtain a divisorial
correspondence on $C\times C$, we choose a point $q\in
\mathbb{P}^1\setminus B$ and denote the fiber $h^{-1}(q)$ by
$\{x,\iota(x)\}$. Write $f^{-1}(x)=\{y_1,y_2,y_3\}$ and
$f^{-1}\bigl(\iota(x)\bigr)=\{z_1,z_2,z_3\}$; then
$p^{-1}(q)=\{y_i+z_j\,|\,i,j=1,2,3\}$ and the identity
$$
D(y_s+z_t)=\sum_{j\neq t}(y_s+z_j)+\sum_{i\neq s}(y_i+z_t)
$$ defines a fixed point free symmetric correspondence $D$
of bidegree $(2,2)$ on $C\times C$. Note that the matrix of
entries $a_{(s,j),(i,t)}$ (for $(s,j),(i,t)\in\{1,2,3\}^2$) given
by
\[ a_{(s,j),(i,t)}=\left\{ \begin{array}{lll}
1 & {\rm if} & (y_s+z_j,\,y_i+z_t)\in D\\
0 & {\rm else}
\end{array}\right. \] is the adjacency matrix of $L_2(3)$.
Hence the canonical endomorphism $\gamma_{\!{}_D}$ of $J(C)$
satisfies the equation
$$
(\gamma_{\!{}_D}-{\rm id}_{J(X)})(\gamma_{\!{}_D}+2\,{\rm
id}_{J(X)})=0\,.
$$ For $l\geq 1$ it follows that
$P:={\rm im}(\gamma_{\!{}_D}-{\rm id}_{J(C)})$ is an
$l$-dimensional Prym-Tyurin variety of exponent 3 for the curve
$C$ of genus $3l+4$. We shall call $P$ (resp.\ $p$) the Prym
variety (resp.\ covering) associated to $f$ and $h$.\\

The preceding construction is a special case of Example
\ref{prymtyurin}. This is a direct consequence of Lemma
\ref{types} and the following result, for integers $l\geq 1$.

\begin{proposition}
\label{comparison} The Lange-Recillas-Rochas family of
$l$-dimensional Prym-Tyurin varieties of exponent 3 coincides with
the family of Prym-Tyurin varieties $P_{\!{}_+}(\mathcal{T})$ for
Prym triples $\mathcal{T}$ of type $l$.
\end{proposition}

\proof Let $P$ (resp.\ $p:C\rightarrow \mathbb{P}^1$) be the Prym
variety (resp.\ covering) associated to an \'{e}tale threefold
covering $f:\tilde{X}\rightarrow X$ and a double covering
$h:X\rightarrow \mathbb{P}^1$ with branch locus $B$ of cardinality
$2l+8$. We fix a point $q\in \mathbb{P}^1\setminus B$ and write
$N=\{1,2,3\}$. Using the notation of the preceding construction,
we define the bijections $\nu:p^{-1}(q)\rightarrow N\times N$,
$y_i+z_j\mapsto (i,j)$ and $\mu:(h\circ f)^{-1}(q)\rightarrow
\{1,2\}\times N$, sending $y_i\mapsto (1,i)$ and $z_j\mapsto
(2,j)$. We let $\rho$ (resp.\ $\varrho$) be the monodromy
representations for $p$ (resp.\ $h\circ f$) induced by $\nu$
(resp.\ $\mu$). Choose a small $q$-based loop $\lambda\subset
\mathbb{P}^1\setminus B$ around a point $b\in B$. Then
$\varrho([\lambda])=\upsilon_1\upsilon_2\upsilon_3$ is the product
of mutually disjoint transpositions
$\upsilon_j=\bigl((1,s_j)\hskip0.25cm (2,t_j)\bigr)$, $s_j,t_j\in
N$. Employing the fact that $p$ comes from $h\circ f\circ\pi$ with
$\pi$ as in the construction, one easily checks that
$\xi\rho([\lambda])\xi^{-1}$ is a conjugate of some
$\sigma_{m,h}\in S_3\times S_3$ (in the notation of Example
\ref{prymtyurin}). By transitivity of ${\rm im}\,\rho$ it thus
follows that $\xi({\rm im}\,\rho)\xi^{-1}=S_3\times S_3$. Hence,
if $A$ denotes the adjacency matrix of the graph $L_2(3)$, then
$\mathcal{T}=(A,p,[\nu])$ is a Prym triple of type $l$ and Lemma
\ref{matrixcorrespondence} implies that
$P=P_{\!{}_+}(\mathcal{T})$.\\

Conversely, let $\mathcal{T}$ be a Prym triple of type $l$
associated to a covering $p:C\rightarrow \mathbb{P}^1$ with branch
locus $B$ and a labelling class $[\nu]$, where
$\nu:p^{-1}(q)\rightarrow \{1,2,3\}^2$ is a labelling for the
fiber of $p$ over a point $q\in \mathbb{P}^1\setminus B$. Given
the monodromy representation $\rho:\pi_1(\mathbb{P}^1\setminus
B,q)\rightarrow {\rm Perm}(\{1,2,3\}^2)$ for $p$ induced by $\nu$,
we take coordinate mappings $\rho_1,\rho_2$ such that
$\rho(\beta)$ splits as $\bigl(\rho_1(\beta),\rho_2(\beta)\bigr)$
for $\beta\in \pi_1(\mathbb{P}^1\setminus B,q)$. We then have a
transitive representation $\varrho:\pi_1(\mathbb{P}^1\setminus
B,q)\rightarrow {\rm Perm}(\{1,2\}\times \{1,2,3\})$, defined by
the relations
$\varrho(\beta)(1,i)=\bigl(2,\rho_2(\beta)(i,1)\bigr)$ and
$\varrho(\beta)(2,j)=\bigl(1,\rho_1(\beta)(1,j)\bigr)$. Using the
local monodromy of $p$, one shows by analogy with the proof of
Theorem \ref{factorization} that $\varrho$ is a monodromy
representation for a covering
$$
h\circ f:\tilde{X}\xrightarrow{3:1} X\xrightarrow{2:1}
\mathbb{P}^1\,,
$$ where $f$ is \'{e}tale and $g(X)\geq 4$. Then $\rho$ is immediately seen
to act as a monodromy representation for the covering that is
associated to $f$ and~$h$. Hence $P_{\!{}_+}(\mathcal{T})$ (resp.\
$p$) is the Prym variety (resp.\ covering) associated to $f$ and
$h$.\qed

\end{section}

\end{document}